\documentclass[english]{article}
\usepackage[T1]{fontenc}
\usepackage[latin9]{inputenc}
\usepackage[letterpaper]{geometry}
\usepackage{array}
\usepackage{textcomp}
\usepackage{amsthm}
\usepackage{amsmath}
\usepackage{amssymb}
\usepackage[authoryear]{natbib}

\makeatletter

\providecommand{\tabularnewline}{\\}

\theoremstyle{plain}
\newcommand{\lyxaddress}[1]{
\par {\raggedright #1
\vspace{1.4em}
\noindent\par}
}
\theoremstyle{plain}
\newtheorem{thm}{Theorem}
  \theoremstyle{definition}
  \newtheorem{defn}[thm]{Definition}
  \theoremstyle{remark}
  \newtheorem{rem}[thm]{Remark}
  \theoremstyle{plain}
  \newtheorem{lem}[thm]{Lemma}
  \theoremstyle{plain}
  \newtheorem{prop}[thm]{Proposition}
  \theoremstyle{plain}
  \newtheorem{cor}[thm]{Corollary}
 \theoremstyle{definition}
  \newtheorem{example}[thm]{Example}

\makeatother

\usepackage{babel}

\begin{document}

\title{A frequentist framework of inductive reasoning}

\author{David R. Bickel}

\maketitle

\lyxaddress{Ottawa Institute of Systems Biology\\
Department of Biochemistry, Microbiology, and Immunology\\
Department of Mathematics and Statistics \\
University of Ottawa\\
451 Smyth Road\\
Ottawa, Ontario, K1H 8M5}

\lyxaddress{+01 (613) 562-5800, ext. 8670\\
dbickel@uottawa.ca}
\begin{abstract}
Reacting against the limitation of statistics to decision procedures,
R. A. Fisher proposed for inductive reasoning the use of the fiducial
distribution, a parameter-space distribution of epistemological probability
transferred directly from limiting relative frequencies rather than
computed according to the Bayes update rule. The proposal is developed
as follows using the confidence measure of a scalar parameter of interest.
(With the restriction to one-dimensional parameter space, a confidence
measure is essentially a fiducial probability distribution free of
complications involving ancillary statistics.)

A betting game establishes a sense in which confidence measures are
the only reliable inferential probability distributions. The equality
between the probabilities encoded in a confidence measure and the
coverage rates of the corresponding confidence intervals ensures that
the measure's rule for assigning confidence levels to hypotheses is
uniquely minimax in the game. 

Although a confidence measure can be computed without any prior distribution,
previous knowledge can be incorporated into confidence-based reasoning.
To adjust a \textit{p}-value or confidence interval for prior information,
the confidence measure from the observed data can be combined with
one or more independent confidence measures representing previous
agent opinion. (The former confidence measure may correspond to a
posterior distribution with frequentist matching of coverage probabilities.)
The representation of subjective knowledge in terms of confidence
measures rather than prior probability distributions preserves approximate
frequentist validity. 
\end{abstract}

\newpage{}

\section{Introduction}

Within the history of frequentism, the preeminent precedent for basing
inductive inferences directly on hypothesis probabilities is that
found in fiducial arguments of the mature Fisher \citep{RefWorks:1215,RefWorks:1179,RefWorks:52,RefWorks:193,FraserLindley},
who objected as staunchly against the behavioristic theory of Neyman
and Pearson as against the Bayesian theories of Laplace and Jeffreys
(\citealp{RefWorks:1326}; \citeyear{RefWorks:985}). Several interpretations
of fiducial inference have been advanced. For example, \citet{RefWorks:305}
attempted to justify it in terms of a derivation of logical parameter
probabilities from statistical coverage rates under a principle of
irrelevance defined in terms of conditional independence \citep{RefWorks:1163}.
Fiducial-like arguments are now employed in the context of functional
models \citep[e.g.,][]{Kohlas2008378} and generalized inference \citep[e.g.,][]{RefWorks:1175}.
The context of the present study is that of an intelligent agent formulating
a probabilistic level of certainty of a hypothesis on the basis of
observed data and possibly on the basis of more subjective information
as well.

A \emph{significance function} is a cumulative distribution function
(CDF) of a probability distribution called a \emph{confidence measure},
which is equivalent to a fiducial probability measure if ancillarity
considerations are neglected. The confidence measure succinctly includes
all the information needed to compute any confidence intervals or
\textit{p}-values of any null hypotheses for a single scalar parameter
of interest. As will be seen, such combination methodology allows
the utilization of one or more subjective probability distributions
of that parameter without forfeiting frequentist validity. As in strict
Bayesian inference, each subjective distribution represents the knowledge
of some expert or other intelligent agent, allowing the principled
incorporation of existing information into an analysis of observed
data. The use of such information in a frequentist framework in the
case of a scalar parameter of interest is made possible by recent
methodology for combining significance functions, originally intended
for the meta-analysis of independent data sets \citep{RefWorks:130}.

\subsection{\label{sub:Significance-and-confidence}Confidence measures}

An observed sample $x\in\Omega$ of $n$ observations is modeled as
a realization of the random quantity $X$ of probability distribution
$P_{\xi},$ where $\xi$ is the value of the parameter of some family
of distributions. Letting $\theta=\theta\left(\xi\right)\in\Theta$
denote the subparameter of interest and $\gamma=\gamma\left(\xi\right)\in\Gamma$
the nuisance subparameter, $\left\langle \theta,\gamma\right\rangle $
is written in place of $\xi$ without loss of generality.

\begin{defn}[Significance function]
\label{def:CD}The function \[
F:\Omega\times\Theta\rightarrow\left[0,1\right]\]
is a \emph{significance function }for $\theta=\theta\left(\xi\right)$
if \[
F\left(x,\bullet\right)=F_{x}\left(\bullet\right):\Theta\rightarrow\left[0,1\right]\]
is a cumulative distribution function (CDF) for all $x\in\Omega$
and if \begin{equation}
P_{\left\langle \theta,\gamma\right\rangle }\left(F_{X}\left(\theta\right)<\alpha\right)=\alpha\label{eq:uniform}\end{equation}
for all $\theta\in\Theta$, $\gamma\in\Gamma$, and $\alpha\in\left[0,1\right]$.
\end{defn}
The condition of equation (\ref{eq:uniform}) says that $F_{X}\left(\theta\right)$
is a pivotal quantity with a uniform distribution on $\left[0,1\right].$ 

\noindent \noindent The significance function encodes a rich set of confidence
intervals, as follows.

\noindent \begin{lem}
\label{lem:CI}If $F$ is a significance function with inverse function
$F^{-1}:\Omega\times\left[0,1\right]\rightarrow\Theta$, then\begin{equation}
P_{\left\langle \theta,\gamma\right\rangle }\left(\theta\in\left(F_{X}^{-1}\left(\alpha_{1}\right),F_{X}^{-1}\left(1-\alpha_{2}\right)\right]\right)=1-\alpha_{1}-\alpha_{2}\label{eq:CD}\end{equation}
for all $\theta\in\Theta,$ $\gamma\in\Gamma,$ and $\alpha_{1},\alpha_{2}\in\left[0,1\right]$
such that $\alpha_{1}+\alpha_{2}\le1$. Conversely, consider the function
$F^{-1}:\Omega\times\left[0,1\right]\rightarrow\Theta$ such that
$F_{x}^{-1}$ is an inverse CDF for all $x\in\Omega$. If equation
(\ref{eq:CD}) holds for all $\theta\in\Theta,$ $\gamma\in\Gamma,$
and $\alpha_{1},\alpha_{2}\in\left[0,1\right]$ such that $\alpha_{1}+\alpha_{2}\le1$,
then $F:\mathcal{B}\times\Omega\rightarrow\left[0,1\right]$, the
inverse of $F^{-1},$ is a significance function.
\end{lem}
The straightforward proof is omitted. The significance function provides
standard one- and two-sided \emph{p}-values for testing the null hypothesis
that $\theta=\theta^{\prime}$ as well as exact $\left(1-\alpha_{1}-\alpha_{2}\right)100\%$
confidence intervals. A test against the alternative hypothesis $\theta>\theta^{\prime},$
$\theta<\theta^{\prime},$ or $\theta\ne\theta^{\prime}$ respectively
yields $F_{X}\left(\theta^{\prime}\right)$, $1-F_{X}\left(\theta^{\prime}\right)$,
or $2F_{X}\left(\theta^{\prime}\right)\wedge2\left(1-F_{X}\left(\theta^{\prime}\right)\right)$
as the \emph{p}-value \citep{RefWorks:60,RefWorks:127}.

The significance function is used to generate a probability measure
of the interest parameter:
\begin{defn}[confidence measure]
\label{def:certainty-distribution}Consider $F$, a significance
function for $\theta$. For all $x\in\Omega,$ if $F_{x}$ is the
CDF of a random quantity $\vartheta$ that has some probability distribution
$P^{x}$ on the measurable space $\left(\Theta,\mathcal{B}\right)$,
then $P^{x}$ is the \emph{confidence measure} of\emph{ }$\theta$
that corresponds to $F$ given $X=x$.
\end{defn}

\citet{Efron19933} dubbed $P^{x}$ a \emph{confidence distribution},
the term \citet{RefWorks:127} and \citet{RefWorks:130} instead attached
to $F$ due to the isomorphism noted below. To avoid confusion between
the probability measure and its CDF, $F$ is herein called the significance
function, following \citet{RefWorks:60}. For clarity, the emphasis
will be on the probability distribution rather than on the significance
function since confidence measures take the place of Bayesian prior
and posterior measures. The idea of a confidence distribution goes
as far back as \citet{conditionalityPrinciple1958}, who recommended
the simultaneous consideration of confidence intervals for multiple
levels of confidence. \citet{Polansky2007b}, referring to $P^{x}$
probabilities as attained or observed confidence levels, provides
an accessible introduction to the concept and its applications. 

The connection between a confidence measure and such confidence intervals
is directly made in this statement that the inferential probability
that $\vartheta$ is in a particular observed confidence interval
is equal to the coverage rate of the random confidence interval that
it realizes.
\begin{lem}
\label{lem:logical-probability-is-CI-rate}Given a random quantity
$\vartheta$ that has some confidence measure $P^{x}$ on $\left(\Theta,\mathcal{B}\right)$
corresponding to $F$ given $X=x$, \begin{eqnarray}
1-\alpha_{1}-\alpha_{2} & = & P^{x}\left(\vartheta\in\left(F_{x}^{-1}\left(\alpha_{1}\right),F_{x}^{-1}\left(1-\alpha_{2}\right)\right]\right)\label{eq:certainty-interval}\\
 & = & P_{\left\langle \theta,\gamma\right\rangle }\left(\theta\in\left(F_{X}^{-1}\left(\alpha_{1}\right),F_{X}^{-1}\left(1-\alpha_{2}\right)\right]\right)\end{eqnarray}
for all $x\in\Omega,$ $\theta\in\Theta,$ $\gamma\in\Gamma,$ and
$\alpha_{1},\alpha_{2}\in\left[0,1\right]$ such that $\alpha_{1}+\alpha_{2}\le1$.\end{lem}
\begin{proof}
Exact frequentist coverage at rate $1-\alpha_{1}-\alpha_{2}$ follows
from Lemma \ref{lem:CI}. That rate is equal to the parameter-space
probability given $X=x$: \begin{eqnarray*}
1-\alpha_{1}-\alpha_{2} & = & F_{x}\left(F_{x}^{-1}\left(1-\alpha_{2}\right)\right)-F_{x}\left(F_{x}^{-1}\left(\alpha_{1}\right)\right)\\
 & = & P^{x}\left(\vartheta\le F_{x}^{-1}\left(1-\alpha_{2}\right)\right)-P^{x}\left(\vartheta<F_{x}^{-1}\left(\alpha_{1}\right)\right).\end{eqnarray*}

\end{proof}
This result will be generalized to arbitrary confidence sets in Section
\ref{sub:Set-estimation}. 

A confidence measure can be constructed from any significance function:
\begin{lem}
\label{lem:extension}Given some significance function $F$, there
is a random quantity $\vartheta$ of a confidence measure $P^{x}$
that corresponds to $F$ given $X=x$ such that, for all $\theta\in\Theta$
and $x\in\Omega,$\begin{equation}
F_{x}\left(\theta\right)=P^{x}\left(\vartheta<\theta\right).\label{eq:extension}\end{equation}
\end{lem}
\begin{proof}
For all $x\in\Omega,$ consider a function $\tilde{P}^{x}:\mathcal{B}\rightarrow\left[0,1\right]$
that satisfies $\tilde{P}^{x}\left(\left(\theta^{\prime},\theta^{\prime\prime}\right]\right)=F_{x}\left(\theta^{\prime\prime}\right)-F_{x}\left(\theta^{\prime}\right)$
for all $\theta^{\prime},\theta^{\prime\prime}\in\Theta$ such that
$\theta^{\prime}\le\theta^{\prime\prime}.$ By the Caratheodory extension
theorem (e.g., \citealp[pp. 578-581]{Schervish1995} or \citealp[pp. 26-27]{Kallenberg1}),
there is a measure space $\left(\Theta,\mathcal{B},P^{x}\right)$
such that $\tilde{P}^{x}\left(\Theta^{\prime}\right)=P^{x}\left(\Theta^{\prime}\right)$
for all $\Theta^{\prime}\in\mathcal{B}$. Then $P^{x}$ is a confidence
measure corresponding to $F$ given $X=x$ with the random quantity
$\vartheta:\Theta\rightarrow\Theta$. 
\end{proof}

Thus, every significance function evaluated at $X=x$ is isomorphic
to a confidence measure.

\subsection{\label{sub:agent}Intelligent agents}

For the sake of clarity and close contact with actual problems of
statistical data analysis, decision-theoretic results will be presented
in familiar terms of estimation rather than solely in terms of abstract
decision makers. It is nonetheless often expedient to refer to such
hypothetical agents to place the present work in context with the
literature since many have found it convenient to imagine an ideally
information-processing agent such as the robot of \citet{RefWorks:1108},
especially when motivating axiomatic decision theory and its game-theoretic
precursors (§\ref{sub:Axiomatic-decision-theory}). While algorithmic
agents in artificial intelligence often make real decisions, agents
in statistics instead inform a researcher or administrator who will
consider the data analysis results and their underlying assumptions
when making a decision that cannot be completely automated. To avoid
confusion with actual people, pronouns and possessives referring to
agents will be neuter.

\subsection{\label{sub:Overview}Overview}

Section \ref{sec:Properties-of-confidence} establishes properties
of confidence measures that both motivate their use and guide their
subjective assignment. Various definitions and lemmas of Sections
\ref{sub:Significance-and-confidence} and \ref{sub:Confidence-based-estimation}
provide a framework for the accounts of coherence in Section \ref{sec:Coherence-and-rationality}
and for a game-theoretic attribute of the confidence measure that
gives precise, general content to the following reasoning. \citet[p. 224]{Kempthorne-Jaynes}
considered fair odds for betting on the hypothesis that an observed
confidence interval covers the parameter value to be a function of
the rate of frequentist coverage $\rho$ as if he were using a confidence
measure $P^{x}$, claiming that such a betting strategy would outperform
a Bayesian, {}``coherently wrong'' strategy. Heuristically, the
thought is that in assessing a fair betting rate, achieving a reported
frequency of correct decisions over repeated sampling outweighs the
importance of coherence over time; cf. \citet{RefWorks:1161}. The
rational component of Kempthorne's assertion had been formally specified
in terms of minimizing risk under a simple loss function \citep{RefWorks:1187}.
That risk is generalized to a risk associated with testing arbitrary
hypotheses in Section \ref{sub:Reliability}, which establishes that
the only minimax solutions are confidence measures. 

Section \ref{sec:Incorporation} turns to the special case of subjective
confidence measures as defined in Section \ref{sec:Subjective-confidence-measure}.
A strategy of combining confidence measures, including one of more
subjective measures, is proposed in Section \ref{sub:Combining}.
Guidance on the assignment of subjective confidence measures is then
given in Section \ref{sub:Assignment-of-subjective}. Subjective confidence
may be assigned to hypotheses (1) indirectly by means of a hypothetical
data set on which the agent might rest its opinion, (2) directly on
the basis of minimaxity or other frequentist properties of the confidence
measure, or (3) indirectly by transforming a Bayesian prior distribution
into a confidence measure.

The paper concludes with a brief discussion.

\section{\label{sec:Properties-of-confidence}Properties of confidence measures}

Sections \ref{sec:Coherence-and-rationality} and \ref{sub:Reliability}
respectively record coherence and decision-theoretic criteria met
by confidence measures using the terminology introduced in Section
\ref{sub:Confidence-based-estimation}.

\subsection{\label{sub:Confidence-based-estimation}Confidence-based estimation}

Section \ref{sub:Set-estimation} treats the problem of deriving a
subset $\hat{\Theta}\left(x\right)$ of the parameter set $\Theta$
that has a desired level of confidence or rate of coverage $\rho$;
the set $\hat{\Theta}\left(x\right)$ is called a \emph{set estimate}
to distinguish it from a point estimate such as that of Section \ref{sub:Hypothesis-indicator-estimation}.

\subsubsection{\label{sub:Hypothesis-indicator-estimation}Hypothesis indicator
estimation}

The degree to which a hypothesis is considered supported by data is
defined as a point estimate of the value indicating whether the hypothesis
is true:
\begin{defn}
\label{def:indicator-estimator}A function \begin{equation}
\hat{1}:\mathcal{B}\times\Omega\rightarrow\left[0,1\right]\label{eq:indicator-estimator}\end{equation}
is called an \emph{indicator estimator} on $\mathcal{B}\times\Omega$.
For all $\theta\in\Theta$, $\Theta^{\prime}\in\mathcal{B}$, and
$x\in\Omega$, the value $\hat{1}\left(\Theta^{\prime},x\right)$,
hereafter written as $\hat{1}_{\Theta^{\prime}}\left(x\right),$ is
an \emph{estimate} of $1_{\Theta^{\prime}}\left(\theta\right).$\end{defn}
\begin{rem}
\label{rem:indicator-estimator}This agrees with the interpretation
of inferential or logical probability as an estimate of the truth
value of its hypothesis \citep[e.g., ][]{Wilkinson1977,RefWorks:1203}.
However, the definition is general enough to include in principle
any function from $\mathcal{B}\times\Omega$ to $\mathbb{R}^{1}$
by use of a monotonic transform to the conventional $\left[0,1\right]$
range. 
\end{rem}
Evaluating the indicator estimator under squared error loss, \citet{Hwang1992490}
found that $F_{\bullet}\left(\theta^{\prime}\right)$ and $1-F_{\bullet}\left(\theta^{\prime}\right)$
are admissible estimators of $1_{\left(\inf\Theta,\theta^{\prime}\right)}\left(\theta\right)$
and $1_{\left(\theta^{\prime},\sup\Theta\right)}\left(\theta\right)$,
respectively, in the case of exponential models, with $\theta$ as
the location parameter. The resulting squared-error admissibility
of $P^{\bullet}\left(\vartheta<\theta^{\prime}\right)$ as an estimator
of $1_{\left(\inf\Theta,\theta^{\prime}\right)}\left(\theta\right)$
is a weak condition satisfied by all generalized Bayes rules \citep{Hwang1992490}
regardless of their actual frequentist performance.

\subsubsection{\label{sub:Set-estimation}Set estimation}

In order to lay the groundwork for the minimax result of Section \ref{sub:Reliability},
a general set estimator is defined in terms of the general indicator
estimator in the same way as confidence intervals are often defined
in terms of \emph{p}-values. Let $\lambda$ denote the Lebesgue measure
on $\mathbb{R}^{1}$ and $\mathcal{B}\left(\left[0,1\right]\right)$
the Borel $\sigma$-field of $\left[0,1\right].$ 
\begin{defn}
\label{def:set-estimator}A function $\hat{\Theta}:\mathcal{B}\left(\left[0,1\right]\right)\times\Omega\rightarrow\mathcal{B}$
is a \emph{set estimator} and, if the map $\hat{\Theta}_{\bullet}\left(x\right):\mathcal{B}\left(\left[0,1\right]\right)\rightarrow\mathcal{B}$
is bijective for all $x\in\Omega$, then $\hat{\Theta}$ is an \emph{invertible}
\emph{set estimator}. Further, $\hat{\Theta}$ is the \emph{set estimator}
corresponding to an indicator estimator $\hat{1}$ on $\mathcal{B}\times\Omega$
if $\hat{\Theta}$ is a set estimator and if $\hat{1}_{\hat{\Theta}_{B}\left(x\right)}\left(x\right)=\lambda\left(B\right)$
for all $x\in\Omega$. Each observed $\hat{\Theta}_{B}\left(x\right)$
is a \emph{set estimate}; $\lambda\left(B\right)$ is the \emph{level}
or \emph{nominal probability} of a \emph{particular set estimator}
$\hat{\Theta}_{B}$ with index $B$ in $\mathcal{B}\left(\left[0,1\right]\right)$.
\end{defn}
The confidence coefficient and Bayesian credibility are examples of
the level $\lambda\left(B\right)$ of a particular set estimator.
Each set $B$ in $\mathcal{B}\left(\left[0,1\right]\right)$ is used
to index a particular set estimator in order to facilitate working
with $\left\{ \hat{\Theta}_{B}:B\in\mathcal{B}\left(\left[0,1\right]\right)\right\} $,
a comprehensive collection of particular set estimators corresponding
to the same indicator estimator. This proves more convenient than
indexing particular set estimators with their levels since the same
level can correspond to multiple particular set estimators. For example,
the lower-tail $\left(B=\left[0,0.95\right)\right)$, upper-tail $\left(B=\left(0.05,1\right]\right)$,
and central $\left(B=\left(0.025,0.975\right)\right)$ 95\% Bayesian
credibility intervals represent three particular set estimators, each
of the same level, 95\%. Since $B$ is a Borel set and since $\lambda$
is the Lebesgue measure on $\mathbb{R}^{1},$ less usual indices such
as $B=\left(0.05,0.10\right)\cup\left(0.50,0.99\right)$ are also
possible.

The following lemma and theorem are also needed for the game-theoretic
result of the next section.
\begin{lem}
\label{lem:rate-is-level}Suppose there are some significance function
$F$ and indicator estimator $\hat{1}$ on $\mathcal{B}\times\Omega$
such that\begin{equation}
\hat{1}_{\left(\theta^{\prime},\theta^{\prime\prime}\right]}\left(x\right)=F_{x}\left(\theta^{\prime\prime}\right)-F_{x}\left(\theta^{\prime}\right),\label{eq:deltaCD}\end{equation}
for all $\theta^{\prime},\theta^{\prime\prime}\in\Theta$ such that
$\theta^{\prime}\le\theta^{\prime\prime}$ and for all $x\in\Omega.$
If \textup{$\hat{\Theta}_{B}:\mathcal{B}\left(\left[0,1\right]\right)\times\Omega\rightarrow\mathcal{B}$}
is an invertible set estimator corresponding to $\hat{1}$, then \begin{equation}
\lambda\left(B\right)=P_{\left\langle \theta,\gamma\right\rangle }\left(\theta\in\hat{\Theta}_{B}\left(X\right)\right)\label{eq:rate-is-level}\end{equation}
for all $B\in\mathcal{B}\left(\left[0,1\right]\right)$, $\theta\in\Theta$
and $\gamma\in\Gamma$. Conversely, if there is an invertible set
estimator $\hat{\Theta}:\mathcal{B}\left(\left[0,1\right]\right)\times\Omega\rightarrow\mathcal{B}$
corresponding to an indicator estimator $\hat{1}$ on $\mathcal{B}\times\Omega$
such that equation (\ref{eq:rate-is-level}) holds for all $B\in\mathcal{B}\left(\left[0,1\right]\right)$,
$\theta\in\Theta,$ and $\gamma\in\Gamma,$ then there is some significance
function $F$ such that $\hat{1}$ satisfies equation (\ref{eq:deltaCD})
for all $\theta^{\prime},\theta^{\prime\prime}\in\Theta$ such that
$\theta^{\prime}\le\theta^{\prime\prime}$ and for all $x\in\Omega.$
\end{lem}
\begin{proof}
According to Lemma \ref{lem:CI}, exact coverage (\ref{eq:rate-is-level})
holds for every set estimator $\hat{\Theta}_{B}\left(X\right)$ that
maps to an interval subset of $\mathcal{B}.$ To prove exact coverage
of every set estimator $\hat{\Theta}_{B}\left(X\right)$ that maps
to a union of disjoint interval subsets of $\mathcal{B}$, note that
$\hat{\Theta}_{B^{\prime}}\left(x\right)$ is the subset of $\hat{\Theta}_{B}\left(x\right)$
corresponding to subset $B^{\prime}$ of $B$ for some $x\in\Omega$
according to the invertibility of $\hat{\Theta}$. With $\mathcal{B}\left(\left[0,1\right]\right)^{\prime}\left(B\right)$
denoting the set of all disjoint interval subsets of $B$, \begin{eqnarray*}
P_{\left\langle \theta,\gamma\right\rangle }\left(\theta\in\hat{\Theta}_{B}\left(X\right)\right) & = & \sum_{B^{\prime}\in\mathcal{B}\left(\left[0,1\right]\right)^{\prime}\left(B\right)}P_{\left\langle \theta,\gamma\right\rangle }\left(\theta\in\hat{\Theta}_{B^{\prime}}\left(X\right)\right)\\
 & = & \sum_{B^{\prime}\in\mathcal{B}\left(\left[0,1\right]\right)^{\prime}\left(B\right)}\lambda\left(B^{\prime}\right)=\lambda\left(B\right)\end{eqnarray*}
for all $B\in\mathcal{B}\left(\left[0,1\right]\right)$, $\theta\in\Theta$
and $\gamma\in\Gamma,$ thereby proving the first half of the lemma.
The converse follows from Lemma \ref{lem:CI} and the fact that equation
(\ref{eq:CD}) is a special case of equation (\ref{eq:rate-is-level}).
\end{proof}
Equation (\ref{eq:rate-is-level}) says the level of any particular
set estimator is equal to the actual coverage rate of that set estimator.
Hence, the probability that $\vartheta$ is in a particular estimated
set is equal to the coverage rate of the corresponding set estimator:
\begin{thm}
\label{thm:logical-probability-is-rate}Suppose there are some indicator
estimator $\hat{1}$ on $\mathcal{B}\times\Omega$ and some significance
function $F$ such that, for all $\Theta^{\prime}\in\mathcal{B}$
and $x\in\Omega,$\begin{equation}
\hat{1}_{\Theta^{\prime}}\left(x\right)=P^{x}\left(\vartheta\in\Theta^{\prime}\right),\label{eq:probability}\end{equation}
where $\vartheta$ is a random quantity of measure $P^{x},$ the confidence
measure of $\theta$ given $X=x$ that corresponds to $F.$ Let $\hat{\Theta}:\mathcal{B}\left(\left[0,1\right]\right)\times\Omega\rightarrow\mathcal{B}$
denote any invertible set estimator corresponding to $\hat{1}.$ Then
\begin{equation}
P^{x}\left(\vartheta\in\hat{\Theta}_{B}\left(x\right)\right)=\lambda\left(B\right)=P_{\left\langle \theta,\gamma\right\rangle }\left(\theta\in\hat{\Theta}_{B}\left(X\right)\right)\label{eq:logical-probability-is-rate}\end{equation}
for all $x\in\Omega,$ $B\in\mathcal{B}\left(\left[0,1\right]\right)$,
$\theta\in\Theta$ and $\gamma\in\Gamma$. Conversely, if there is
an invertible set estimator $\hat{\Theta}:\mathcal{B}\left(\left[0,1\right]\right)\times\Omega\rightarrow\mathcal{B}$
corresponding to an indicator estimator $\hat{1}$ on $\mathcal{B}\times\Omega$
such that \textup{$\lambda\left(B\right)=P_{\left\langle \theta,\gamma\right\rangle }\left(\theta\in\hat{\Theta}_{B}\left(X\right)\right)$
}for all $B\in\mathcal{B}\left(\left[0,1\right]\right)$, $\theta\in\Theta,$
and $\gamma\in\Gamma,$ then there is some significance function $F$
and some confidence measure of $\theta$ given $X=x$ that corresponds
to $F$ such that equations (\ref{eq:probability}) and (\ref{eq:logical-probability-is-rate})
hold for all $x\in\Omega,$ $B\in\mathcal{B}\left(\left[0,1\right]\right)$,
$\theta\in\Theta$ and $\gamma\in\Gamma$.\end{thm}
\begin{proof}
By Lemma \ref{lem:logical-probability-is-CI-rate}, equation (\ref{eq:logical-probability-is-rate})
holds for all interval elements of $\mathcal{B}\left(\left[0,1\right]\right).$
That result for intervals is extended to all unions of disjoint intervals
in $\mathcal{B}\left(\left[0,1\right]\right)$ by the invertibility
as used in the proof of Lemma \ref{lem:rate-is-level}, thereby proving
the first half of the theorem. The converse follows directly from
Lemmas \ref{lem:rate-is-level} and \ref{lem:extension}.
\end{proof}

Succinctly generalizing equation (\ref{eq:logical-probability-is-rate})
to $\mathbf{\theta}$ as a vector parameter of interest, \citet[pp. 4-5, 69, 224-227]{Polansky2007b}
defined rather than derived $P^{x}\left(\Theta^{\prime}\right),$
the {}``attained confidence level'' of $\theta\in\Theta^{\prime},$
to be the coverage frequency of a corresponding confidence set $\hat{\Theta}_{\rho,\omega}\left(X\right):$
\begin{equation}
P^{x}\left(\Theta^{\prime}\right)=\rho=P_{\left\langle \theta,\gamma\right\rangle }\left(\theta\in\hat{\Theta}_{\rho,\omega\left(\rho\right)}\left(X\right)\right),\label{eq:weight-is-rate}\end{equation}
where the coverage rate $\rho$ and shape parameter $\omega\left(\rho\right)$
are constrained such that $\hat{\Theta}_{\rho,\omega\left(\rho\right)}\left(x\right)=\Theta^{\prime}$
for the observed value $x$ of random element $X,$ the distribution
$P_{\left\langle \theta,\gamma\right\rangle }$ of which is indexed
by parameter $\left(\left\langle \theta,\gamma\right\rangle \right).$
The {}``observed confidence level'' \citep{Polansky2007b} of equation
(\ref{eq:weight-is-rate}) should not be confused with theories of
estimating confidence levels \citep{RefWorks:1277,GoutisCasella1995b}.

\subsection{\label{sec:Coherence-and-rationality}Axiomatic coherence}

Each of the next two subsections establishes the coherence of the
certainty distribution $P^{x}$ from a distinct viewpoint that led
to axioms of coherence or rationality. The first perspective is decision-theoretic
and the second is logic-theoretic.

\subsubsection{\label{sub:Axiomatic-decision-theory}Axiomatic decision theory}

\paragraph{\label{sub:Dutch}Precursors to axiomatic decision theory}

The use of the certainty distribution for decision making was motivated
in Section \ref{sub:Reliability} by placing the decision-making agent
in the role of a casino that will settle bets at its published betting
odds, allowing a gambling opponent to choose hypotheses on which to
bet. This represents situations in which an agent must make a definite
decision on the basis of limited information, as when it must either
accept the hypothesis that the true parameter value is a pre-specified
interval or accept the hypothesis that is is in the complement of
that interval. 

That is essentially the gambling scenario for which Ramsay and de
Finetti considered this \emph{Dutch book} situation: a gambler can
contract bets with any casino agent that assesses betting odds for
certain events in violation of probability theory such that the agent
will lose regardless of the outcomes \citep[pp. 59-65]{RefWorks:Gillies2000}.
An agent or indicator estimator $\hat{1}$ is called \emph{coherent}
if it assigns betting odds in such a way that it will not suffer such
\emph{sure loss}. \citet{Schervish1995} presents the equivalent mathematical
definition of coherence to which the following proposition refers.
\begin{prop}
\label{lem:coherent}Let $\mathcal{M}$ be the collection of all measurable
maps from a measurable space $\left(\Omega,\Sigma\right)$ to $\left(\Theta,\mathcal{B}\right)$.
An indicator estimator $\hat{1}$ on $\mathcal{B}\times\Xi$ is coherent
if and only if there is a probability measure $P$ on $\left(\Omega,\Sigma\right)$
such that $\hat{1}_{\Theta^{\prime}}\left(x\right)=P\left(\vartheta\in\Theta^{\prime}\right)$
for all $\vartheta\in\mathcal{M}$ and $\Theta^{\prime}\in\mathcal{B}$.\end{prop}
\begin{proof}
This follows immediately from \citet[Theorem B.139]{Schervish1995},
who uses notation and terminology closer to that of \citet{RefWorks:43}.
\end{proof}
The definition of conditional probability has been recovered by a
similar theorem based on bets that are called off if some event does
not occur; see, e.g., \citet[pp. 657-658]{Schervish1995} or \citet{Hacking2001}.
In an idealized framework, setting conditional betting rates by any
parameter distribution other than a conditional probability distribution
leads to certain loss \citep{RefWorks:1168,RefWorks:1187,RefWorks:1185,RefWorks:1167,RefWorks:1165}.
Since the probabilities in the theorems provided have no time dependence,
they do not indicate the method of replacing a parameter distribution
after new data are observed and thus are compatible with the proposed
method of replacement by maintaining correct confidence interval coverage
rates (§\ref{sec:frequentist-framework}). In Bayesian inference,
on the other hand, the parameter distribution used to place bets after
observing data is identified with the prior distribution conditional
on the observed data. Such identification is an assumption that is
usually hidden, not a consequence of coherence \citep{CoherentFrequentism}.

There are known problems with resting coherence on Dutch book theorems
alone \citep{RefWorks:1253,Howson2009177}. De Finetti admitted that
arguments from betting behavior do not provide an unobjectionable
foundation for coherent decision making \citep{RefWorks:Gillies2000}.
Ramsay also looked beyond the Dutch book argument, speculating that
an axiomatic foundation encompassing both utility and probability
could be laid \citep[p. 30]{French2000}. \citet{RefWorks:126} proved
the conjecture by drawing on the game theory conceived in mathematical
economics, and others have since created generalizations of his axiomatic
decision theory \citep{French2000}.

\paragraph{Axiomatic decision theory proper}

Although axiomatic systems of decision theory were developed with
subjective probability in mind, nothing in the mathematics prohibits
more objective applications by interpreting hypothesis probabilities
as indicator estimates rather than as levels of belief. In fact, the
axioms only put very weak constraints on rational decision-making
that lead to coherently representing unknown values as random quantities
without requiring the additional constraints of a prior distribution
and the characteristically Bayesian use of conditional probability.
In place of the latter constraints, the proposed framework substitutes
the requirement that probabilities correspond to frequentist rates
of coverage. 

While specifying a particular utility function for use with the axioms
is inherently subjective, it is no more so than specifying a particular
loss function for use in classical frequentist decision theory or
a particular significance or confidence level for use in Neyman-Pearson
theory. In order to objectively communicate the results of data analysis,
probability distributions of parameters can be reported without utilities,
as is common Bayesian practice. Accordingly, reporting a certainty
distribution of a parameter allows each agent to supply its own loss
function when making decisions on the basis of what can be inferred
about the parameter value from the available data.

\subsubsection{\label{sub:Axiomatic-inductive-logic}Axiomatic inductive logic}

While the axiomatic decision theories, building on foundations laid
by Bayes \citep[§1.3]{RefWorks:182}, Ramsay, and de Finetti, derive
probability from the maximization of expected utility rather than
vice versa (§\ref{sub:Axiomatic-decision-theory}), many have questioned
the propriety of the order \citep[e.g., ][]{RefWorks:1214}. That
order was reversed by \citet{RefWorks:1127}, \citet{RefWorks:182},
Cox \citeyearpar{Cox19461,RTCox}, \citet{RefWorks:612}, and \citet{RefWorks:1201},
who constructed axiomatic formulations of inductive-logical probability
on parameter space without relying on betting behavior, expected gain,
or other decision-theoretic concepts. 

The term \emph{logical probability} is used here in the broad sense
of mathematical probability interpreted according to any axiomatic
system that generalizes some logic of deduction. Because such systems
have been closely associated with some version of the now discredited
principle of insufficient reason \citep[p. 64]{RefWorks:999,RefWorks:Gillies2000},
the statistical community has not deemed them a practical guide for
data analysis. Logical probability may prove more useful in practice
when supplemented instead with a frequentist principle such as one
of minimizing arbitrary-hypothesis risk (\ref{eq:risk}).

The system of Cox \citeyearpar{Cox19461,RTCox}, while lacking mathematical
rigor, remains highly regarded its simplicity and for the generality
of its assumptions \citep[e.g., ][]{Paris232410,RefWorks:999,VanHorn20033,Howson2009177}
and continues to convince scientists to express uncertainty probabilistically
\citep[e.g., ][]{PhysRevE.72.031912}. Its two axioms may be expressed
in the notation of Section \ref{sub:Reliability} with the addition
of joint and conditional indicator estimators $\hat{1}$ in the second
axiom \citep[pp. 3-4]{RTCox}:
\begin{enumerate}
\item $\hat{1}_{\Theta\backslash\Theta^{\prime}}\left(x\right)$ is a smooth
function of $\hat{1}_{\Theta^{\prime}}\left(x\right)$ for all $x\in\Xi$
and $\Theta^{\prime}\subseteq\Theta$.
\item $\hat{1}_{\Theta^{\prime},\Theta^{\prime\prime}}\left(x\right)$,
the estimate of $1_{\Theta^{\prime}}\left(\theta\right)\wedge1_{\Theta^{\prime\prime}}\left(\theta\right),$
is a smooth function of $\hat{1}_{\Theta^{\prime}}\left(x\right)$
and of $\hat{1}_{\Theta^{\prime\prime}}\left(x|\theta\in\Theta^{\prime}\right)$,
the conditional estimate of $1_{\Theta^{\prime\prime}}\left(\theta\right)$
given $\theta\in\Theta^{\prime}$, for all $x\in\Xi,$ $\emptyset\subset\Theta^{\prime}\subseteq\Theta,$
and $\Theta^{\prime\prime}\subseteq\Theta$.
\end{enumerate}
From more general versions of those stated axioms, a few tacit assumptions,
and the rules of classical logic, \citet{RTCox} proved $\hat{1}$
to be isomorphic to finitely additive probability \citep{Paris232410,VanHorn20033,Howson2009177},
allowing identification with the certainty distribution (\ref{eq:probability-yields-correct-odds})
as well as with the Bayesian posterior that Cox originally had in
mind.

\subsection{\label{sub:Reliability}Game-theoretic interpretation}

\subsubsection{Decisions as bets on hypotheses\label{sub:non-Bayesian-posterior}}

Caution is needed when drawing general conclusions from the losses
suffered by gambling agents since such conclusions can be sensitive
to the rules of the game \citep{RefWorks:1160}. Further, some games
resemble situations faced in practice better than others. By construction,
inference according to the proposed methodology is robust across two
games so different that each had been used to argue for an opposite
paradigm of statistics:
\begin{enumerate}
\item \citet{Kempthorne-Jaynes} and \citet{RefWorks:1285} alluded to a
game like that of Section \ref{sub:Reliability} to support Neyman-Pearson
statistics;
\item The game of posting fair betting odds for and against every hypothesis
in some $\sigma$-field is the foundation of the traditional Dutch-book
argument for Bayesian statistics. See \citet[Appendix C]{Purdue}
for an accessible summary of \citet[pp. 654-655]{Schervish1995},
who furnishes a more general theorem. 
\end{enumerate}

\subsubsection{\label{sec:frequentist-framework}Arbitrary-hypothesis minimaxity}

While pure Neyman-Pearson inference is optimal under a risk function
that in effect imposes an infinite penalty for failing to control
a Type I error rate at some specified level, such a risk function
does not provide a helpful representation of all situations faced
by the statistician. Many situations that call for data-based decisions
are better represented by a risk function representing a statistician's
necessity to give odds for the hypothesis that an observed confidence
interval covers the parameter of interest such that a decision maker
can use those odds to safely bet either for or against that hypothesis
as directed by an opponent \citep{RefWorks:1187}; this game gives
structure to the claims of \citet{Kempthorne-Jaynes} and \citet{RefWorks:1285}
that were mentioned in Section \ref{sub:Overview}. 

That risk function is extended to accommodate more general hypothesis
testing via the following zero-sum game played between a statistician
and a client. The client will specify a pair of mutually exclusive
and jointly exhaustive hypotheses to which the statistician must assign
betting odds. Those odds determine the amount of either a payoff or
penalty for the statistician, depending on which hypothesis is true. 

This situation is further stylized by representing the decision-making
statistician as a casino agent and the opposing client as a gambler
at the casino. The statistician applies a comprehensive collection
of set estimators to data and, for each level-$\rho$\emph{ }set estimate,
posts $\rho/\left(1-\rho\right)$ as fair betting odds for the event
that the set estimate includes $\theta$ to the event that the set
estimate does not include $\theta$. The set is comprehensive in the
sense that its elements map to all elements of $\mathcal{B}$ for
each $x\in\Omega$. In posting fair betting odds, the statistician
announces a willingness commit to paying the client $\rho$ or less
if the set estimate does not include $\theta$ provided that $1-\rho$
or more would instead be received from the client if the set estimate
includes $\theta$. The statistician also must swap the payment amounts
to bet that the set estimate does not include $\theta$ if the client
desires. The client only accepts bet proposals at the odds the statistician
considers fair, not favorable. Further, knowing the distributions
of the set estimators in the statistician's set, the client will not
accept unfavorable bets, that is, bets with negative risk to the statistician.
The client enforces this by computing $\omega$, the truly fair betting
odds as defined by the ratio of the rate at which sets from the statistician
cover $\theta$ to the rate of its non-coverage. The client then compares
the fair betting odds to $\rho/\left(1-\rho\right)$ when deciding
whether to accept a bet at odds $\rho/\left(1-\rho\right)$. Thus,
the statistician only successfully contracts a bet on coverage if
$\rho/\left(1-\rho\right)\ge\omega$ or on non-coverage if $\rho/\left(1-\rho\right)\le\omega$.
This contract is concisely represented in terms of loss suffered by
the statistician:\[
\mathcal{L}_{B}\left(\hat{\Theta};X\right)=\begin{cases}
\rho1_{\Theta\backslash\hat{\Theta}_{B}\left(X\right)}\left(\theta\right)-\left(1-\rho\right)1_{\hat{\Theta}_{B}\left(X\right)}\left(\theta\right), & \rho/\left(1-\rho\right)>\omega\left(B\right)\\
\left(1-\rho\right)1_{\hat{\Theta}_{B}\left(X\right)}\left(\theta\right)-\rho1_{\Theta\backslash\hat{\Theta}_{B}\left(X\right)}\left(\theta\right), & \rho/\left(1-\rho\right)<\omega\left(B\right)\\
0, & \rho/\left(1-\rho\right)=\omega\left(B\right),\end{cases}\]
where $B\in\mathcal{B}\left(\left[0,1\right]\right),$ $\rho=\lambda\left(B\right),$
and $\hat{\Theta}$ is an invertible set estimator mapping $\mathcal{B}\left(\left[0,1\right]\right)\times\Omega$
to $\mathcal{B}$ and corresponding to some indicator estimator $\hat{1}$
on $\mathcal{B}\times\Omega$; the \emph{fair betting odds} of $\theta\in\hat{\Theta}_{B}\left(X\right)$
to $\theta\notin\hat{\Theta}_{B}\left(X\right)$ are given by \begin{equation}
\omega\left(B\right)=\omega_{\left\langle \theta,\gamma\right\rangle }\left(B\right)=\frac{P_{\left\langle \theta,\gamma\right\rangle }\left(\theta\in\hat{\Theta}_{B}\left(X\right)\right)}{P_{\left\langle \theta,\gamma\right\rangle }\left(\theta\notin\hat{\Theta}_{B}\left(X\right)\right)},\label{eq:odds}\end{equation}
resulting in the risk the statistician assumes by relying on $\hat{1}$
for assessing the odds of an arbitrary hypothesis.
\begin{defn}
\label{def:risk}Consider $\mathcal{D}\left(\hat{1}\right)$, the
collection of all invertible set estimators each mapping $\mathcal{B}\left(\left[0,1\right]\right)\times\Omega$
to $\mathcal{B}$ and corresponding to some indicator estimator $\hat{1}$
on $\mathcal{B}\times\Omega$. The \emph{arbitrary-hypothesis risk
}of $\hat{1}$ is \begin{eqnarray}
R_{\left\langle \theta,\gamma\right\rangle }\left(\hat{1}\right) & =\min_{\hat{\Theta}\in\mathcal{D}\left(\hat{1}\right)}\max_{B\in\mathcal{B}\left(\left[0,1\right]\right)} & E_{\left\langle \theta,\gamma\right\rangle }\left(\mathcal{L}_{B}\left(\hat{\Theta};X\right)\right)\label{eq:risk}\end{eqnarray}
for all $\theta\in\Theta$ and $\gamma\in\Gamma.$
\end{defn}
As in Neyman-Pearson testing, the hypotheses to be assessed are arbitrary
in the sense that they are dictated by the needs of the current application
and are thus outside of the agent's control. Additional arbitrary
hypotheses may also be specified in the future for unforeseen applications.
For the purpose of defining the risk associated with the indicator
estimator $\hat{1}$ used to assess an arbitrary hypothesis, the worst-case
specification of a hypothesis corresponds to the least-favorable selection
of the corresponding set estimator $\hat{\Theta}_{B}$. Derivation
of a testing procedure from a set estimator rather than vice versa
is not without precedent \citep{RefWorks:1260,Liu1997266,RefWorks:249,RefWorks:1137}.
\begin{lem}
\label{lem:exact-coverage}The indicator estimator $\hat{1}$ on $\mathcal{B}\times\Omega$
is minimax to arbitrary-hypothesis risk if and only if there is an
invertible set estimator $\hat{\Theta}:\mathcal{B}\left(\left[0,1\right]\right)\times\Omega\rightarrow\mathcal{B}$
corresponding to $\hat{1}$ such that\begin{equation}
P_{\left\langle \theta,\gamma\right\rangle }\left(\theta\in\hat{\Theta}_{B}\left(X\right)\right)=\lambda\left(B\right)\label{eq:correct coverage}\end{equation}
 for all $B\in\mathcal{B}\left(\left[0,1\right]\right),$ $\theta\in\Theta$,
and $\gamma\in\Gamma.$ \end{lem}
\begin{proof}
An indicator estimator $\hat{1}$ is minimax to arbitrary-hypothesis
risk if and only if it minimizes $\max_{\theta\in\Theta,\gamma\in\Gamma}R_{\left\langle \theta,\gamma\right\rangle }\left(\hat{1}\right)$
(Definition \ref{def:risk}). Given a particular set estimator $\hat{\Theta}_{B}$
for any $B\in\mathcal{B}\left(\left[0,1\right]\right),$ the odds
for $\theta\in\hat{\Theta}_{B}$ are $\lambda\left(B\right)/\left(1-\lambda\left(B\right)\right)$
as assessed by $\hat{1}$. If equation (\ref{eq:correct coverage})
holds, those odds are equal to $\omega\left(B\right)$, the true odds
given by equation (\ref{eq:odds}), and thus $E_{\left\langle \theta,\gamma\right\rangle }\left(\mathcal{L}_{B}\left(\hat{\Theta};X\right)\right)=0$
for all $B\in\mathcal{B}\left(\left[0,1\right]\right),$ $\theta\in\Theta,$
and $\gamma\in\Gamma.$ Therefore, $\max_{\theta\in\Theta,\gamma\in\Gamma}R_{\left\langle \theta,\gamma\right\rangle }\left(\hat{1}\right)=0$.
But if equation (\ref{eq:correct coverage}) does not hold, then $\lambda\left(B\right)/\left(1-\lambda\left(B\right)\right)\ne\omega\left(B\right)$.
If $\exists B\in\mathcal{B}\left(\left[0,1\right]\right),\theta\in\Theta,\gamma\in\Gamma$
such that $\lambda\left(B\right)/\left(1-\lambda\left(B\right)\right)>\omega\left(B\right),$
then \begin{eqnarray*}
E_{\left\langle \theta,\gamma\right\rangle }\left(\mathcal{L}_{B}\left(\hat{\Theta};X\right)\right) & = & \lambda\left(B\right)P_{\left\langle \theta,\gamma\right\rangle }\left(\theta\notin\hat{\Theta}_{B}\left(X\right)\right)-\left(1-\lambda\left(B\right)\right)P_{\left\langle \theta,\gamma\right\rangle }\left(\theta\in\hat{\Theta}_{B}\left(X\right)\right)\end{eqnarray*}
\[
\frac{E_{\left\langle \theta,\gamma\right\rangle }\left(\mathcal{L}_{B}\left(\hat{\Theta};X\right)\right)}{\left(1-\lambda\left(B\right)\right)P_{\left\langle \theta,\gamma\right\rangle }\left(\theta\notin\hat{\Theta}_{B}\left(X\right)\right)}=\frac{\lambda\left(B\right)}{1-\lambda\left(B\right)}-\omega\left(B\right)>0.\]
Likewise, if $\exists B\in\mathcal{B}\left(\left[0,1\right]\right),\theta\in\Theta,\gamma\in\Gamma$
such that $\lambda\left(B\right)/\left(1-\lambda\left(B\right)\right)<\omega\left(B\right),$
then\[
\frac{E_{\left\langle \theta,\gamma\right\rangle }\left(\mathcal{L}_{B}\left(\hat{\Theta};X\right)\right)}{\left(1-\lambda\left(B\right)\right)P_{\left\langle \theta,\gamma\right\rangle }\left(\theta\notin\hat{\Theta}_{B}\left(X\right)\right)}=\omega\left(B\right)-\frac{\lambda\left(B\right)}{1-\lambda\left(B\right)}>0\]
for all $\theta\in\Theta,\gamma\in\Gamma$. Both results together
indicate that if there is any $B$ in $\mathcal{B}\left(\left[0,1\right]\right)$
such that equation (\ref{eq:correct coverage}) does not hold for
any $\theta\in\Theta,\gamma\in\Gamma$, then $\max_{\theta\in\Theta,\gamma\in\Gamma}R_{\left\langle \theta,\gamma\right\rangle }\left(\hat{1}\right)>0$. 
\end{proof}
That lemma leads to the corollary of Theorem \ref{thm:logical-probability-is-rate}
that establishes the unique minimaxity of the confidence measure.
\begin{cor}
\label{cor:probability-yields-correct-odds}The indicator estimator
$\hat{1}$ on $\mathcal{B}\times\Omega$ is minimax to arbitrary-hypothesis
risk if and only if there is some significance function $F$ such
that, for all $\Theta^{\prime}\in\mathcal{B}$ and $x\in\Omega,$\begin{equation}
\hat{1}_{\Theta^{\prime}}\left(x\right)=P^{x}\left(\vartheta\in\Theta^{\prime}\right),\label{eq:probability-yields-correct-odds}\end{equation}
where $\vartheta$ is a random quantity of law $P^{x},$ the confidence
measure of $\theta$ that corresponds to $F$ given $X=x$.\end{cor}
\begin{proof}
By Lemma \ref{lem:exact-coverage}, this corollary obtains if and
only if exact coverage (\ref{eq:correct coverage}) holds for every
particular set estimator $\hat{\Theta}_{B}$ corresponding to the
indicator estimator $\hat{1}$ given by equation (\ref{eq:probability-yields-correct-odds}).
Theorem \ref{thm:logical-probability-is-rate} supplies the necessary
and sufficient conditions. 

 \end{proof}
\begin{rem}
The conditions of Theorem \ref{thm:logical-probability-is-rate} and
Corollary \ref{cor:probability-yields-correct-odds} include continuous
data and exact satisfaction of the uniformity condition for brevity
and clarity. Applications to actual data often require judicious approximations.
For example, a {}``half-correction'' makes the significance function
applicable to discrete data \citep{RefWorks:127}. As an alternative
to approximate confidence measures, upper and lower probabilities
constituting envelopes of a class of confidence measures have been
proposed \citep{CoherentFrequentism} in the spirit of the Dempster-Shafer
theory of belief functions \citep{RefWorks:1221}. In that case, methods
of combining confidence measures (§\ref{sub:Combining}; \citealp{RefWorks:130})
would apply separately to each confidence measure of the class.
\end{rem}

\section{\label{sec:Incorporation}Incorporation of previous information}

\subsection{\label{sec:Subjective-confidence-measure}Subjective confidence measures}

Consider a set $\mathcal{A}$ of agents, a set $\mathcal{T}$ of
$\Omega$-measurable maps, and a function $t_{\bullet}:\mathcal{A}\rightarrow\mathcal{T}.$
The random quantity $t_{a}\left(X\right)$ represents relevant information
accessible to any agent $a\in\mathcal{A}$ or what that agent can
recall from previous experience. The information is \emph{subjective
}in the sense that \[
\begin{aligned}t_{a^{\prime}}\left(x\right)\ne t_{a^{\prime\prime}}\left(x\right), &  & a^{\prime},a^{\prime\prime}\in\mathcal{A}\end{aligned}
\]
is possible for the same observed data $x\in\Omega.$ Errors in agent
perception or memory could be modeled by, for each $a\in\mathcal{A},$
modeling $t_{a}$ as the realization of a random function $T_{a},$
e.g., yielding \[
T_{a}\left(X\right)=\left\langle X_{1}+Y\left(a\right),X_{2}+Y\left(a\right),\dots,X_{n}+Y\left(a\right)\right\rangle ,\]
where $Y\left(a\right)$ is some random variable independent of $X_{i}.$
Likewise, agents themselves may be randomly selected from $\mathcal{A},$
in which case $A$ represents a random agent. 

For generality, $\mathbf{P}_{\theta,\gamma}$ will represent a joint
probability measure extended from $P_{\left\langle \theta,\gamma\right\rangle }$
such that $\mathbf{P}_{\theta,\gamma}\left(X=\bullet\right)$ is the
distribution of the data on which each agent indirectly relies, $\mathbf{P}_{\theta,\gamma}\left(A=\bullet\right)$
is the distribution of agents, and $\mathbf{P}_{\theta,\gamma}\left(T=\bullet\right)$
is the distribution of maps from $\Omega$ to $\mathcal{T}.$ The
trivial assignments $\mathbf{P}_{\theta,\gamma}\left(A=a\right)=1$
and $\mathbf{P}_{\theta,\gamma}\left(T=t\right)=1$ are important
special cases. If there is a function $G:\mathcal{T}\times\Theta\rightarrow\left[0,1\right]$
such that \begin{equation}
\mathbf{P}_{\theta,\gamma}\left(G_{T_{A}\left(X\right)}\left(\theta\right)<\alpha\right)=\alpha\label{eq:subjective-measure}\end{equation}
for all $\theta\in\Theta$, $\gamma\in\Gamma$, and $\alpha\in\left[0,1\right],$
then $G$ is a significance function by Definition \ref{def:CD}.
Let $Q^{\left\langle x,a,t\right\rangle }$ denote the confidence
measure of $\theta$ that corresponds to $G$ given $\left\langle X,A,T\right\rangle =\left\langle x,a,t\right\rangle $.
Because $G$ and $Q^{\left\langle x,a,t\right\rangle }$ depend on
subjective information in the form of $T_{A},$ they are called a
\emph{subjective significance function }and a \emph{subjective confidence
measure}, respectively. By contrast, a confidence measure that does
not depend on such subjective information is called \emph{objective}.

\subsection{\label{sub:Combining}Updating a subjective confidence measure}

\subsubsection{Combining subjective and objective measures}

\noindent Subjective Bayesianism makes use of agent knowledge by deriving
the prior distribution of the parameters conditional on the data.
The result is the Bayes posterior distribution of the parameters,
from which a marginal posterior distribution of the parameter of interest
may be obtained. The proposed substitute for this application of Bayes's
theorem is the \textit{combination} of the subjective confidence measure
with the objective confidence measure by generating a new confidence
measure. Due to the isomorphism between confidence measures and significance
functions (§\ref{sub:Significance-and-confidence}), any of the methods
of combining significance functions studied by \citet{RefWorks:130}
is equivalent to combining their corresponding confidence measures.
Thus, I propose that such methods be used to combine subjective and
objective confidence measures by combining their significance functions
into a single significance function corresponding to a combined confidence
measure.
\begin{example}
\citet{RefWorks:130} proposed a method of significance function combination
that relies on the choice of a continuous cumulative distribution
function. Out of several such functions considered, $DE,$ the cumulative
distribution function of the double exponential distribution, performs
the best under repeated sampling, especially for small sample sizes
\citep{RefWorks:130}, when the consideration of agent opinion has
the most impact. For this function, the significance function $F,$
and $L$ independent samples $X\left(1\right),X\left(2\right),...,X\left(L\right)$
each drawn from $P_{\left\langle \theta,\gamma\right\rangle },$ the
combined significance function $\tilde{F}$ is defined such that

\begin{equation}
\tilde{F}\left(\theta\right)=DE_{L}\left(DE^{-1}\left(F_{X\left(1\right)}\left(\theta\right)\right)+DE^{-1}\left(F_{X\left(2\right)}\left(\theta\right)\right)+\cdots+DE^{-1}\left(F_{X\left(L\right)}\left(\theta\right)\right)\right)\label{ZEqnNum134425}\end{equation}
for all $\theta\in\Theta,$ $DE_{L}\left(q\right)$ is the convolution
of \textit{L} copies of $DE\left(q\right),$ and

\[
DE^{-1}\left(p\right)=\log\left(2p\right)I_{}\]
with \textit{I} as the indicator function. \citet{RefWorks:130} found
that the convolution may be computed from $V_{L}\left(q\right),$
a polynomial that satisfies a simple recursive relation, using

\[
DE_{L}\left(q\right)=\left(1-\frac{1}{2e^{q}}V_{L}\left(q\right)\right)I_{\left[0,\infty\right)}\left(q\right)+\frac{1}{2e^{-q}}V_{L}\left(-q\right)I_{\left(-\infty,0\right]}\left(q\right).\]
The most commonly used polynomials are $V_{1}\left(q\right)=1$, $V_{2}\left(q\right)=1+\frac{q}{2}$,
and $V_{3}\left(q\right)=1+\frac{5q+q^{2}}{8}$.
\end{example}
A confidence measure formed by combining both objective and subjective
confidence measures will be called an \textit{agent-updated confidence
measure}. Such formation will be referred to as \emph{agent-based
confidence updating} (ABCU).

\subsubsection{Reducing sensitivity to violations of matching}

\noindent Since the distribution of $T_{a}$ is typically unknown,
an available subjective distribution is at best an approximation to
a confidence measure. Combining subjective distributions with objective
confidence measures can be made more robust to the influence of a
subjective distribution's deviation from the properties of a confidence
measure by treating a subjective distribution as if it were an objective
confidence measure based on an incorrect model. In effect, a confidence
measure derived from an incorrect model or from a misleading agent
opinion is not a confidence measure of the true parameter value, but
of some other underlying parameter value. Methods of adaptive confidence-measure
combination assign weights such that a combination of confidence measures
with different underlying parameter values asymptotically has the
properties of a confidence measure with respect to the true parameter
value \citep{RefWorks:130}. Thus the value of the asymptotic combined
confidence measure, when evaluated at the true value of the parameter,
has the uniform distribution (\ref{eq:uniform}) needed.

\subsubsection{\label{sub:common-mean}Multiple sampling distributions}

Confidence-measure combination applies independent samples drawn from
different populations as well as to independent samples drawn from
the same population. The following example illustrates this.
\begin{example}
\label{exa:common-mean}A standard common-mean problem involves estimating
a mean $\theta$ shared by two or more normal populations of unknown
variances that may differ from one population to another. Confidence
measure combination using equation \eqref{ZEqnNum134425} yielded
95\% confidence intervals with close to 95\% coverage for simulated
data drawn from two normal populations of $\theta=1$ for various
choices of sample sizes and population variances \citep{RefWorks:130}.
For one such choice, samples of sizes 3 and 4 were drawn from populations
of variance 1 and $3.5^{2},$ respectively. For this illustration,
those populations generated the realizations $y_{1}=\left(0.523,2.460,1.119\right)$
and $y_{2}=\left(0.072,-2.275,-4.554,-0.077\right)$, two observed
samples with means $\left(\bar{x}_{1},\bar{x}_{2}\right)$ and standard
deviations $\left(s_{1},s_{2}\right)$. Their combined confidence
measure $\tilde{F}$ then gives

\[
\tilde{F}\left(\theta\right)=DE_{2}\]
according to equation \eqref{ZEqnNum134425}. Since the sample sizes
are small, any agent knowledge may play an influential role in inference
about the mean. Given that the subjective distributions on $\theta$
for two independent agents are ${\rm N}\left(0,3^{2}\right)$ and
${\rm N}\left(2,4^{2}\right)$, the agent-updated confidence measure
gives

\[
\tilde{F}_{a}\left(\theta\right)=DE_{4},\]
where $\Phi\left(q\right)$ is the standard normal cumulative distribution
function. $\tilde{F}_{a}\left(\theta\right)$ is a confidence measure
under the assumption of a matching agent space. The effect of including
agent opinion may be quantified by noting, e.g., the impact on the
\textit{p}-value of $H_{0}:\theta_{{\rm true}}=-1$ versus $H_{a}:\theta_{{\rm true}}>-1$;
$\tilde{F}\left(-1\right)\approx0.104$, whereas $\tilde{F}_{a}\left(-1\right)\approx0.049$.
Like in the case of suitably assigned prior probabilities, the effect
of incorporating agent knowledge disappears as either sample size
goes to infinity. 

Writing from a Bayesian perspective, \citet{Seidenfeld2007b} challenged
frequentism with a common-mean problem in the guise of the inferring
the mass of a hollow cube from a direct measurement of the mass of
the cube and a measurement of the mass of a ruler of known density.
In what is arguably the most mature logical system of direct inference,
the consideration of both measurements leads to a loss of precision
compared to consideration of either measurement alone \citep{Kyburg2007b}.
By contrast, confidence-measure combination enables the effective
use not only of information in both measurements but also of any subjective
information, as in Example \ref{exa:common-mean}. 
\end{example}

\subsection{\label{sub:Assignment-of-subjective}Assignment of subjective confidence}

This section outlines three general methods of assigning a subjective
confidence measure. 

However, it is not always necessary to obtain an entire confidence
measure. For example, one subjective \textit{p}-value suffices to
obtain the agent-based \textit{p}-value with respect to a single null
hypothesis.

\subsubsection{Hypothetical-data confidence assignment}

One way to assign a Bayesian prior is to have the agent generate a
hypothetical data set on which its knowledge about $\theta$ is based;
see, e.g., \citet{LeleElicitation}. Likewise, using equation (\ref{eq:subjective-measure}),
a subjective confidence measure $Q^{\left\langle x,a,t\right\rangle }$
can be derived from the hypothetical data $x$ on which agent $a$
might have based its opinion. For that calculation, $t$ either may
be an identity function such that $t\left(x\right)=x$ and $\mathbf{P}_{\theta,\gamma}\left(T=t\right)=1$
or may calibrate the data assignment to correct an elicitation bias.
The technique will be illustrated with a simple example of generating
hypothetical data to generate a subjective likelihood ratio. 
\begin{example}
\label{exa:Edwards}\citet{RefWorks:52}, making use of the likelihood
principle without relying on prior distributions, suggested the use
of subjective likelihood ratios under the paradigm in which inference
is made directly from the actual likelihood function rather than from
a posterior distribution derived from it or from likelihood functions
based on unobserved data. His Example 3.5.1 supposes the knowledge
of Torricelli regarding the atmospheric pressure $\mu$, prior to
his taking measurements, is equivalent to the knowledge that would
be gained by drawing 740 mmHg from a normal distribution of unknown
mean $\mu$ and known standard deviation 25 mmHg. Then, up to a proportionality
constant, the subjective likelihood is $\exp\left(-\frac{1}{2}\left(\frac{740-\mu}{25}\right)^{2}\right)$,
whereas the likelihood from his measurement, 760 mmHg, is $\exp\left(-\frac{1}{2}\left(\frac{760-\mu}{1}\right)^{2}\right)$
since the sampling distribution is assumed normal with standard deviation
1 mmHg. Multiplying the two likelihood functions yields the likelihood
function used for inference; its maximum occurs at $\mu=759.968$
mmHg \citep{RefWorks:52}. 
\end{example}
The same method can generate a subjective confidence measure. 
\begin{example}
In the problem of Example \ref{exa:Edwards}, the subjective and objective
likelihood functions are proportional to the densities of the subjective
and objective confidence measures if $\mathbf{P}_{\theta,\gamma}\left(A=a\right)=1$
and $\mathbf{P}_{\theta,\gamma}\left(T=t\right)=1,$ where $t$ is
an identity function. Using equation (\eqref{ZEqnNum134425}), the
density of the agent-updated confidence measure reaches its maximum
at $\mu=759.231$ mmHg.
\end{example}

\subsubsection{\label{sub:Direct-confidence-assignment}Direct confidence assignment}

Like Bruno de Finetti's prior probability of a hypothesis, the subjective
confidence of a hypothesis can be defined as the perceived value of
the opportunity to gain one unit of utility if the hypothesis is true.
In other words, both the Bayesian subjective probability and the subjective
confidence level are betting quotients that determine the decisions
of some agent. 

While Bayesian prior probability is constrained only by coherence,
confidence levels are also constrained by the rules of the game described
in Section \ref{sub:Reliability}, modified as follows. In reporting
set estimators to the client casino, the statistician agent $a$ cannot
use the data $X$ but rather reduced data $t_{a}\left(X\right)$ defined
by some $\Omega$-measurable map $t_{a}.$ A set estimator $\hat{\Theta}_{B}$
is only available to the agent if there is some function $\Theta^{\dagger}$
such that \begin{equation}
\hat{\Theta}_{B}\left(x\right)=\hat{\Theta}_{B}^{\dagger}\left(t_{a}\left(x\right)\right)\label{eq:reduction}\end{equation}
for all $B\in\mathcal{B}\left(\left[0,1\right]\right),x\in\Omega.$
The agent will achieve minimaxity under arbitrary-hypothesis loss
if and only if the set estimates made on the basis of its choices
of $\hat{\Theta}_{B}^{\dagger}$s are isomorphic to confidence measures.
With the set estimators selected by the agent under restriction (\ref{eq:reduction}),
each confidence measure $Q^{\left\langle x,a,t\right\rangle }$ is
that agent's subjective confidence measure on the basis of information
$t\left(x\right),$ where $t=T$ if $T$ is a fixed function or $t$
is a realization of $T$ if $T$ is a random function. 

As in Bayesian statistics, subjective distributions need not be constructed
by a formal process of eliciting the actual beliefs of a human individual
or organization. Nonetheless, the definition of subjective confidence
in terms of an agent's betting rate serves as a guiding principle
for the specification of subjective confidence measures. 

The foundational idea of this frequentist methodology may be motivated
by the corresponding Bayesian methodology. From an idealized Bayesian
standpoint, inasmuch as an agent's elicited opinion has been coherently
formulated from observed data such as summaries of data more directly
observed by others, the agent's prior will equal the posterior distribution
that would have been obtained by applying Bayes's theorem to the former
data were they still available. Similarly, since under frequentism
a set of \textit{p}-values or confidence intervals takes the place
of the posterior distribution as the inference result, one may consider
eliciting \textit{p}-values or confidence intervals from an agent
that, in the ideal case, equal that would have been computed given
the now unavailable data on which the agent opinion was based. Then,
in analogy with how the Bayesian combines a subjective prior distribution
with a likelihood function obtained from data, the frequentist may
use confidence-measure combination to combine subjective \textit{p}-values
and confidence intervals with those obtained from data. The confidence
measure provides a framework for such combination by encapsulating
the information of \textit{p}-values and confidence intervals into
a confidence measure. Informally, given a sample of fixed observations,
the cumulative distribution function of the standard confidence measure
maps each null hypothesis parameter value to its upper-tailed \textit{p}-value.
Likewise, given a fixed elicitation of an agent's opinion, the cumulative
distribution function of the subjective confidence measure maps each
null hypothesis parameter value to the upper-tailed \textit{p}-value
assigned by the agent. Alternatively, the agent's subjective distribution
may be approximated by interpolating the lower-tailed \textit{p}-values
provided at a sufficient number of null hypothesis parameter values
or by interpolating the confidence intervals provided at a sufficient
number of confidence levels. To the extent that each randomly selected
agent's knowledge is correctly summarized in such a subjective distribution,
one-sided \textit{p}-values derived from those priors follow a uniform
distribution under the truth of the null hypothesis.

\subsubsection{\label{sub:Indirect-confidence-assignment}Bayesian confidence assignment}

Due to the novelty of assigning subjective confidence directly, the
guidance provided in Section \ref{sub:Direct-confidence-assignment}
might be less reliable in the case of a human agent than the application
of mature procedures of eliciting Bayesian prior distributions and
of correcting them to ensure coherence. For examples, see \citet{RefWorks:34},
\citet{RefWorks:39}, and \citet{RefWorks:63}.

If a Bayesian prior elicited from an agent equals a Bayesian posterior
computed from a probability-matching initial distribution \citep{RefWorks:40}
and the data on which the agent indirectly based its prior, then the
prior asymptotically approaches a confidence measure; see \citet{RefWorks:130}
on asymptotic significance functions. In this sense, a Bayesian prior
is approximately equal to a subjective confidence measure suitable
for combination with an objective confidence measure.

\section{Discussion}

\subsection{Reliability of confidence measures}

The reliability of $P^{x}$ in the form of confidence matching (Theorem
\ref{thm:logical-probability-is-rate}) and arbitrary-hypothesis minimaxity
(Corollary \ref{cor:probability-yields-correct-odds}) means the decision-making
agent cannot suffer any expected loss due to a strategy of placing
bets over repeated samples in the above game-theoretic framework.
Consequently, if the game were modified such that at least some of
the bets placed are favorable, the agent would accrue an expected
gain. However, the benefit of achieving minimax risk and the stated
rates of interval coverage is not limited to those rare or non-existent
situations in which more than one sample is drawn from the same population;
rather, those properties reflect the reliability of methods satisfying
them.

\subsection{Other approaches to subjective information}

\subsubsection{Objective and subjective Bayes}

\noindent Considering ABCU as \textit{subjective frequentism} elucidates
its relationship to standard Bayesian approaches (Table \ref{tab:comparison}).
The comparison involves some blurring between Bayesianism and frequentism,
as posterior distributions with probability-matching priors are approximate
confidence measures. This suggests that the choice of whether or not
to incorporate agent knowledge may often tend to have more impact
on inference than the choice of whether to perform Bayesian or frequentist
calculations. 

\begin{table}
\begin{tabular}{|p{1.3in}|p{2in}|p{2in}|}
\hline 
  & \textbf{Objective input only}  & \textbf{Subjective input included} \tabularnewline
\hline 
\textbf{Frequentist calculations}  & Objective confidence measures or combinations thereof are used for
inference  & Combinations of objective and subjective confidence measures are used
for inference \tabularnewline
\hline 
\textbf{Bayesian calculations}  & Starting with an improper prior distribution, the resulting posterior
distribution is used for inference  & Starting with an agent's proper prior distribution, the resulting
posterior distribution is used for inference \tabularnewline
\hline
\end{tabular}\caption{Comparison of subjective Bayesianism and the proposed subjective frequentism
(the agent-updated confidence measure) to their counterparts that
do not formally rely on agent knowledge.\label{tab:comparison}}

\end{table}
Whereas application of the agent-updated confidence measure is only
based on previous information about a one-dimensional parameter of
interest, a subjective Bayesian analysis can also make use of any
information that is also available for the nuisance parameters, often
leading to more reliable inference. However, even when agents do have
such information, it is only rarely elicited, and, in practice, improper
priors tend to be put on the vast majority of nuisance parameters
\citep{RefWorks:15}. That the agent-updated confidence measure requires
the elicitation of information on only one scalar parameter may enable
researchers to incorporate at least some subjective knowledge into
their analyses.

\subsubsection{Previous non-Bayesian approaches}

ABCU is not the only non-Bayesian framework available for use of subjective
information about the scalar parameter of interest. Example \ref{exa:Edwards}
illustrates such use in the direct-likelihood framework further developed
by \citet{RefWorks:122}. The likelihood principle of that framework
is not followed by Neyman-Pearson uses of subjective likelihood, also
called a likelihood penalty \citep{RefWorks:127}.

Perhaps those previous non-Bayesian methods of incorporating agent
knowledge are seldom used because they require elicitation of likelihood
ratios, which, as \citet{RefWorks:123} conceded, are understood by
few scientists. If so, then methods requiring only the elicitation
of either a \textit{p}-value or of a few confidence intervals, as
guided by one of the methods of Section \ref{sub:Assignment-of-subjective},
may have wider appeal.

Also obviating the elicitation of a likelihood function, \citet{RefWorks:61}
consider assigning a subjective prior to a scalar parameter of interest
in the location parameterization. To the extent that its treatment
of nuisance parameters approximates the assignment of diffuse priors
to such parameters, this approach will give results similar to those
of the more classical Bayesian approach. However, the reparameterization
approach is better understood on its own terms: applying Bayes's theorem
to a uniform prior in the location parameterization yields a posterior
that matches frequentist confidence intervals \citep{RefWorks:61},
i.e., the Bayesian posterior is a confidence measure. That other priors
in general bring departures from this matching property distinguishes
this method from ABCU. This does not necessarily indicate the superiority
of the latter, as the preservation of matching probability even under
the use of subjective information comes in exchange for inference
additivity, an issue addressed below.

More generally, one may multiply the probability density function
of a matching prior by a likelihood or pseudo-likelihood that is only
a function of the scalar interest parameter, yielding, after normalization,
a posterior distribution for inference. Following \citet{RefWorks:127},
such a likelihood or pseudo-likelihood function is considered \textit{reduced}
since, unlike the full likelihood function, it does not depend on
the nuisance parameters. If this reduced likelihood is proportional
to the density function of the objective confidence measure, then
the posterior probability will equal the normalized product of the
density of the subjective confidence measure and the density of the
objective confidence measure. (\citet{RefWorks:127} found that the
proportionality property holds for some reduced likelihoods.) Under
certain conditions, that normalized product asymptotically approaches
the combined confidence measure found by equation (\eqref{ZEqnNum134425})
with the substitution of the matching prior distribution for a confidence
measure. In general, however, directly using that equation leads to
more exact confidence measures than those obtained by multiplying
confidence measure densities \citep{RefWorks:130}.

The principle of \textit{inference additivity} mentioned above means
the totality of inferences resulting from several analyses, each based
on different information, is identical to the inference from the single
analysis based on simultaneous consideration of all of the information;
cf. \citet{RefWorks:52}. The loss of inference additivity is the
most obvious drawback of ABCU compared to uses of Bayes's theorem.
In the context of agent-updated confidence measures, the way in which
three or more confidence measures are combined may affect the result,
e.g., combining the combination of two confidence measures with a
third yields a confidence measure that is not necessarily equal to
that of the simultaneous combination of all three. (The combination
of the objective confidence measures of Section \ref{sub:common-mean}
with the combination of both subjective confidence measures results
in a \textit{p}-value of 0.054 instead of the simultaneous-combination
\textit{p}-value of 0.049.) 

It has also been noted that inference based on confidence measures,
unlike that based on Bayesian posterior distributions, violates the
likelihood principle \citep{RefWorks:127}. ABCU shares this violation
with other methods designed to have correct coverage when agent opinion
is not incorporated.

It may be concluded that the optimality of one method over another
will depend largely on the availability of information from agents
and on the relative desirability of each of the inference principles
and frequentist properties. Both methods using generalized confidence
measures share a new way to formalize agent knowledge. The unique
benefit of ABCU is its production of probability statements that are
correct in the frequentist sense even after incorporating that knowledge.

\section{Acknowledgments}

Matthias Kohl kindly provided R code (R Development Core Team 2004)
used to compute the convolution of the double exponential distribution
({}``R-help'' list message posted on 12/23/05). I thank Mark Cooper
for helpful feedback and Jean Peccoud, Mark Whitsitt, Chris Martin,
and Bob Merrill for their support of the seed of this paper at Pioneer
Hi-Bred, International \citep{RefWorks:24}.

\bibliographystyle{elsarticle-harv}
\bibliography{refman}

\end{document}